\documentclass[11pt,oneside]{article}
\title{A continuous path of singular masas in the hyperfinite \IIi factor}
\author{Allan Sinclair\thanks{School of Mathematics, University of Edinburgh, U.K. email: \texttt{a.sinclair@ed.ac.uk}}\  \and Stuart White\thanks{Department of Mathematics, University of Glasgow, U.K. email: \texttt{s.white@maths.gla.ac.uk}}}
\date{}

\usepackage{amsmath}
\usepackage{amssymb}
\usepackage{braket} 
\usepackage{enumerate} 
\usepackage{amsthm}\numberwithin{equation}{section}
\usepackage{graphicx}
\usepackage[T1]{fontenc}
\input{xy}  
\xyoption{all} 
\CompileMatrices 

\theoremstyle{plain} 
\newtheorem{thm}{Theorem}[section]
\newtheorem{lem}[thm]{Lemma}
\newtheorem{cor}[thm]{Corollary}
\newtheorem{prop}[thm]{Proposition}

\theoremstyle{definition} 
\newtheorem{dfn}[thm]{Definition}
\newtheorem{notn}[thm]{Notation}
\newtheorem{eg}[thm]{Example}
\newtheorem{con}[thm]{Construction}

\newcommand{\nm}[1]{\left\|#1\right\|}
\newcommand{\md}[1]{\left\arrowvert#1\right\arrowvert}

\newcommand{\nmit}[1]{\nm{#1}_{\infty,2}}
\newcommand{\floor}[1]{\left\lfloor{#1}\right\rfloor}

\newcommand{\IIi}{${\rm{II}}_1$ }
\newcommand{\tr}{\textrm{tr}}

\newcommand{\R}{R}
\newcommand{\N}{N}
\newcommand{\M}{M}

\newcommand{\A}{A}
\newcommand{\B}{B}
\newcommand{\D}{D}

\newcommand{\J}{J} 

\newcommand{\BH}[1]{\mathbb B\left({#1}\right)} 

\newcommand{\dit}{d_{\infty,2}}

\newcommand{\Norm}[1]{\mathcal N\left(#1\right)}

\newcommand{\ce}[2]{\mathbb E_{#1}\left({#2}\right)} 
\newcommand{\Ce}[1]{\mathbb E_{#1}} 
\newcommand{\CE}[1]{e_{#1}} 

\newcommand{\Proj}{\mathcal P_{\textrm{min}}} 

\newcommand{\puk}[1]{\textrm{Puk}\left({#1}\right)}

\newcommand{\Bb}{\mathcal B} 
\newcommand{\Aa}{\mathcal A} 

\renewcommand{\and}{\textrm{ and }}
\renewcommand{\sp}{\textrm{ }}

\begin{document}
\maketitle

\begin{abstract}
Using methods of R.J.Tauer \cite{Tauer.Masa} we exhibit an
uncountable family of singular masas in the hyperfinite \IIi factor
$\R$ all with Puk\'anszky invariant $\{1\}$, no pair of which are
conjugate by an automorphism of $\R$. This is done by introducing an
invariant $\Gamma(\A)$ for a masa $\A$ in a \IIi factor $\N$ as the
maximal size of a projection $e\in\A$ for which $\A e$ contains
non-trivial centralising sequences for $e\N e$. The masas produced
give rise to a continuous map from the interval $[0,1]$ into the
singular masas in $\R$ equiped with the $\dit$-metric.

A result is also given showing that the Puk\'anszky invariant
\cite{Pukanszky.Invariant} is $\dit$-upper semi-continuous.  As a
consequence, the sets of masas with Puk\'anszky invariant $\{n\}$
are all closed.
\end{abstract}

\section{Introduction}\label{Path.Intro}
The study of maximal abelian self-adjoint von Neumann subalgebras
(masas) in \IIi factors dates back to J.Dixmier \cite{Dixmier.Masa} in
1954, who classified them using normalisers.  Given a masa $\A$ in a
\IIi factor $\N$, the normaliser group $\Norm{\A}$ consists of all
the unitaries $u\in\N$ with $u\A u^*=\A$. The masa A is
\emph{Cartan} if this normaliser group generates $\N$ as a von
Neumann algebra whereas at the other end of the spectrum, $\A$ is
called \emph{singular} if $\Norm{\A}\subset\A$.

Given two Cartan masas $\A$ and $\B$ in the hyperfinite \IIi factor $\R$, there is an automorphism $\theta$ of $\R$ with $\theta(\A)=\B$ (\cite{Connes.CFW}). We say that masas $\A$ and $\B$ with this last property are conjugate via an automorphism of $\R$.  The most sucessful
invariant for distinguishing between non-conjugate singular masas is
that of L.Puk\'anszky \cite{Pukanszky.Invariant}, which he used to
give countably many pairwise non-conjugate singular masas in $\R$.
More recently, E.St{\o}rmer and S.Neshveyev \cite{Stormer.Puk} have used
the Puk\'anszky invariant to exhibit uncountably many pairwise
non-conjugate singular masas in $\R$ and they also give two
non-conjugate singular masas in $\R$ with the same Puk\'anskzy
invariant.  One of our objectives here is to produce uncountably main non-conjugate singular masas in the hyperfinite \IIi factor with the same Puk\'anszky invariant.  This result, stated formally as Theorem
\ref{Path.Con.Cor} below, follows directly from Theorem
\ref{Path.Con.Result}.  

\begin{thm}\label{Path.Con.Cor}
There exist uncountably many singular masas in the hyperfinite \IIi
factor $\R$, each with Puk\'anszky invariant $\{1\}$, such that no
pair of these masas is conjugate by an automorphism of $\R$.
\end{thm}

To show the non-conjugacy of pairs of
masas we look for non-trivial centralising sequences for $\R$
lying in these masas --- the idea used by St{\o}rmer and Neshveyev
in \cite{Stormer.Puk} to distinguish between two singular masas
with Puk\'anszky invariant $\{1\}$.  The presence of non-trivial
centralising sequences inside masas has also been used by A.Connes
and V.Jones \cite{Connes.TwoCartan} to give a factor containing
two non-conjugate Cartan masas, and by V.Jones and S.Popa
\cite{Jones.PropertiesMasas} in the context of non-conjugate
semi-regular masas whose normalisers generate the same
irreducible subfactor of $\R$.

There is a natural metric, $\dit$, on the space of all masas of a
\IIi factor, \cite{Sinclair.PertSubalg}.  The uncountably many
masas we shall produce for Theorem \ref{Path.Con.Cor}, will
actually give us a continuous map from the unit interval, $[0,1]$
into this metric space --- a continuous path of pairwise
non-conjugate singular masas.

In the next section we state some background, defining the metric
$\dit$, the Puk\'anszky invariant and Tauer masas.  In section
\ref{Path.Cont} we discuss the behaviour of the Puk\'anskzy
invariant on limits of sequences of masas, showing that it is upper
semicontinous and that the sets of masas with invariant $\{n\}$ are
all closed (Theorem \ref{Path.Cont.Result}, Corollary
\ref{Path.Cont.Closed}).  Next, in section \ref{Path.Gamma}, we
define a $\Gamma$-invariant for masas using centralising sequences
and establish some basic properties for later use.  It is this
invariant we use in section \ref{Path.Con} to show the non-conjugacy
of the masas we construct to establish Theorem
\ref{Path.Con.Result}, the main result of the paper.  The work in
this paper forms part of sections 3.1 and 3.3 of the second authors
PhD thesis \cite{Saw.Thesis}.

\section{Preliminaries}\label{Path.Prelim}
Let $\N$ be a \IIi factor.  Write $\tr$ for the faithful normal
trace on $\N$, and let $\nm{x}_2=\tr(x^*x)^{1/2}$ be the Hilbert
space norm induced on $\N$ by $\tr$.  Write $L^2(\N)$ for the
completion of $\N$ in this norm. Given a linear map
$\Phi:\N_1\rightarrow\N_2$ between two \IIi factors write
$\nmit{\Phi}$ for the norm of $\Phi$ regarded as a map from $\N_1$
into $L^2(\N_2)$ \cite{Sinclair.StrongSing}, that is
$$
\nmit{\Phi}=\sup\Set{\nm{\Phi(x)}_2|x\in\N_1,\nm{x}\leq1}.
$$
Given a von Neumann subalgebra $\M$ of $\N$, let $\Ce{\M}$ be the
unique trace-preserving normal conditional expectation from $\N$
onto $\M$.  This conditional expectation is obtained by restricting
to $\N$ the orthogonal projection $\CE{\M}$ from $L^2(\N)$ onto
$L^2(\M)$.  In \cite{Sinclair.PertSubalg} a metric, $\dit$, is
introduced on the set of all von Neumann subalgebras of $\N$, by
$$
\dit(\M_1,\M_2)=\nmit{\Ce{\M_1}-\Ce{\M_2}}.
$$
This metric is equivalent to an older metric of E.Christensen
defined in \cite{Christensen.Subalgebras}.  As a consequence the set
of all von Neumann subalgebras equipped with $\dit$ is a complete
metric space, and the subsets of all masas, all singular masas, all
subfactors and all irreduicble subfactors are closed,
\cite{Christensen.Subalgebras}.

To define the Puk\'anskzy invariant \cite{Pukanszky.Invariant} of a
masa in the separable \IIi factor $\N$, we form the standard
representation of $\N$ acting by left multplication on $L^2(\N)$.
Let $\J$ denote the modular conjugation operator on $L^2(\N)$ given
by extending $x\mapsto x^*$ from $\N$.  For each $x\in\N$, $\J x\J$
is the operator of right mutiplication by $x^*$ and $x\mapsto\J
x\J$ is a conjugate linear anti-isomorphism of $\N$ onto $\N'$.
Given a masa $\A$ in $\N$, let $\Aa=(\A\cup\J\A\J)''$ --- an
abelian von Neumann subalgebra of $\BH{L^2(\N)}$, so that $\Aa'$ is
type $\rm{I}$. The orthogonal projection $\CE{\A}$ from $L^2(\N)$
onto $L^2(\A)$ lies in $\Aa$ and $\Aa'\CE{\A}=\Aa\CE{\A}=\A\CE{\A}$ ---
an abelian algebra. The Puk\'anszky invariant is obtained by taking
the type decomposition of $\Aa'(1-\CE{\A})$.  More formally,
$\puk{\A}$ is the subset of $\mathbb N\cup\{\infty\}$ consisting of
all those $n$ for which there is a non-zero projection $p\leq
1-\CE{\A}$ in $\Aa$ such that $\Aa'p$ is type ${\rm{I}}_n$
\cite{Pukanszky.Invariant}.

We shall use the methods of R.J.Tauer \cite{Tauer.Masa} to
construct masas in the hyperfinite \IIi factor $\R$.  The second
author introduced the concept of a Tauer masa in $\R$ in
\cite{Saw.Tauer,Saw.Thesis}.  A masa $\A$ in $\R$ is said to be a
\emph{Tauer masa} if there exists an increasing chain
$(\N_n)_{n=1}^\infty$ of matrix algebras with
$(\bigcup_{n=1}^\infty\N_n)''=\R$, such that $\A\cap\N_n$ is a masa
in $\N_n$ for each $n$.  In this case we write $\A_n$ for
$\A\cap\N_n$ and say for emphasis that $\A$ is \emph{Tauer with
respect to $(\N_n)_{n=1}^\infty$}.  Tauer masas have Puk\'anszky
invariant $\{1\}$, \cite[Theorem 4.1]{Saw.Tauer}. Chains
$(\N_n)_{n=1}^\infty$ of matrix algebras in $\R$ can always be
realised as a tensor products.  More formally, there are finite
dimensional subfactors $(\M_m)_{m=1}^\infty$ of $\R$ such that we
have $\N_n=\bigotimes_{m=1}^n\M_m$, for each $n$.  We use the
notation of \cite{Saw.Tauer,Saw.Thesis} to consider the inclusions
$\A_{n_1}\subset\A_{n_2}$ of approximates of a Tauer masa $\A$ with
respect to the chain $(\N_n)_{n=1}^\infty$.  Let $\Proj(\A_{n_1})$
denote the set of minimal projections of $\A_{n_1}$.  The finite
dimensional approximation $\A_{n_2}$ can then be written as
\begin{equation}\label{Path.Intro.Inclusion}
\A_{n_2}=\bigoplus_{e\in\Proj(\A_{n_1})}e\otimes\A_{n_2,n_1}^{(e)},
\end{equation}
for some masas $\A_{n_2,n_1}^{(e)}$ in $\bigotimes_{m=n_1+1}^{n_2}\M_m$.

In \cite[Theorem 3.2]{Saw.Tauer} a technical criteron was given
for a Tauer masa to be singular in terms of these
$\A_{n_2,n_1}^{(e)}$. We use part of this calculation, which is
essentially Proposition 3.5 of \cite{Saw.Tauer}; the exact
statement given can be found as Proposition 2.2.2 of
\cite{Saw.Thesis}.
\begin{prop}\label{Path.Intro.Sing}
Let $\A$ be a Tauer masa in $\R$ with respect to the subfactors
$(\N_n)_{n=1}^\infty$. If for infinitely many $n_1\in\mathbb N$,
each minimal projection $e$ of $\A_{n_1}$ and $\epsilon>0$, there
is an $n_2>n_1$ and a unitary $w_e\in\A_{n_2,n_1}^{(e)}$ with
$$
\nm{\ce{\A_{n_2,n_1}^{(f)}}{w_e}}_2\leq\epsilon,
$$
for every minimal projection $f\neq e$ in $\A_{n_1}$, then $\A$ is singular.
\end{prop}

\section{Semi-continuity of the Puk\'anszky invariant}\label{Path.Cont}
The key tool in determining the limiting behaviour of the
Puk\'anskzy invariant on sequences of masas is a perturbation
theorem for subalgebras of a \IIi factor \cite[Theorem
6.5]{Sinclair.PertSubalg}, which we state below for the
convenience of the reader.
\begin{thm}[{\cite[Theorem 6.5 (ii)]{Sinclair.PertSubalg}}]\label{Path.Cont.Pert}
If $\A$ and $\B$ are masas in a separable \IIi factor $\N$ with
$\dit(\A,\B)\leq\epsilon$, then there are projections $p\in\A$ and
$q\in\B$, and a unitary $u\in\N$ satisfying
\begin{itemize}
 \item $u(\B q)u^*=\A p$;
 \item $\nm{u-\ce{\B}{u}}_2\leq45\epsilon$;
 \item $\tr(p)=\tr(q)\geq1-(15\epsilon)^2$.
\end{itemize}
\end{thm}

\begin{thm}\label{Path.Cont.Result}
Let $\A_n$ be a sequence of masas in a separable \IIi factor $\N$
converging in the $\dit$-metric to a von Neumann subalgebra $\B$
of $\N$.  This $\B$ is a masa in $\N$, and
\begin{equation}\label{Puk.Def.Cont.1}
\puk{\B}\subset\bigcup_{r=1}^\infty\bigcap_{n=r}^\infty\puk{\A_n}.
\end{equation}
\begin{proof}
That the set of masas is $\dit$-closed is due to E.Christensen in
\cite{Christensen.Subalgebras}.  For each $n$, we apply Theorem
\ref{Path.Cont.Pert} to the pair $(\A_n,\B)$ to obtain projections
$p_n\in\A_n$, $q_n\in\B$ and a unitary $u_n\in\N$ satisfying the
conditions of the theorem.  Take $\B_n=u_n^*\A_n u_n$ --- a masa in
$\N$ which has $\B_n q_n=\B q_n$, by the first property of Theorem \ref{Path.Cont.Pert}.

As $\A_n$ converges to $\B$ in $\dit$, the last property of
Theorem \ref{Path.Cont.Pert} ensures that
$$
\lim_{n\rightarrow\infty}\nm{1-q_n}_2=0.
$$
For any $x\in\N$,
\begin{align*}
\nm{q_nJq_nJx-x}_2=\nm{q_nxq_n-x}_2&\leq\nm{q_nx-x}_2+\nm{q_n(xq_n-x)}_2\\
&\leq\nm{q_n-1}_2(\nm{x}+\nm{q_nx})\\
&\leq 2\nm{x}\nm{q_n-1}_2,
\end{align*}
so that the projections $q_nJq_nJ$ in $\Bb_n\cap\Bb$ converge
strongly to $1$, by density of $\N$ in $L^2(\N)$.

Given some $m\in\puk{\B}$, there must be a central projection
$f\in\Bb=\Bb'\cap\Bb$ with $f\leq 1-\CE{\B}$, such that $\Bb'f$ is
of type ${\rm{I}}_m$.  As $q_nJq_nJf$ converges strongly to $f$ we
must have $q_nJq_nJf\neq 0$ for sufficiently large $n$, those
with $n\geq n_1$ say.  Now
$$
\Bb_n'q_nJq_nJ=\Bb'q_nJq_nJ,
$$
a type $\rm{I}$ von Neumann algebra with centre $\Bb_nq_nJq_nJ=\Bb
q_nJq_nJ$.  For $n\geq n_1$, $q_nJq_nJf$ is a non-zero projection in this centre,
and $\Bb'_nq_nJq_nJf$ is then a central cutdown of $\Bb'f$, so a type
${\rm{I}}_m$ von Neumann algebra.

Observe that $q_n$ and $Jq_nJ$ commute with both $\CE{\B}$ and
$\CE{\B_n}$, as $q_n\in\B\cap\B_n$.  We also have
$q_n\CE{\B_n}=q_n\CE{\B}$ and $Jq_nJ\CE{\B_n}=Jq_nJ\CE{\B}$, as
$\B_nq_n=\B q_n$.  In this way, $q_nJq_nJf\leq 1-\CE{\B_n}$, so that
$m\in\puk{\B_n}$, for $n\geq n_1$. As $\B_n$ and $\A_n$ are
unitarily equivalent, $m\in\puk{\A_n}$ for all $n\geq n_1$, exactly
as required.
\end{proof}
\end{thm}

In the special case when the Puk\'anszky invariant of each $\A_n$ is
$\{n\}$, the only possibility for the Puk\'anszky invariant of the
limit masa $\B$ is also $\{n\}$.
\begin{cor}\label{Path.Cont.Closed}
Let $\N$ be a separable \IIi factor. For each $n\in\mathbb
N\cup\{\infty\}$, the set of all masas with Puk\'anszky invariant
$\{n\}$ is $\dit$-closed.
\end{cor}

In general we do not have equality in ({\ref{Puk.Def.Cont.1}).
\begin{eg}\label{Path.Cont.Example}
Let $\A$ be a masa in the hyperfinite \IIi factor $\R$ with
Puk\'anskzy invariant $\{1\}$.  Take projections $p_n\neq 1$ in $\A$ with
$p_n\rightarrow 1$ strongly.  For
each $n$, let $\B_n$ be a masa in the hyperfinite \IIi factor
$(1-p_n)\R(1-p_n)$ with Puk\'anszky invariant $\{2\}$.  The
existance of such masas dates back to Puk\'anszky's original
examples in \cite{Pukanszky.Invariant}.  Define
$$
\A_n=\Set{ap_n+b|a\in\A,b\in\B_n},
$$
which is a masa in $\R$.  It is then immediate that
$\dit(\A_n,\A)\rightarrow 0$ as $n\rightarrow\infty$ and that both
$1$ and $2$ lie in $\puk{\A_n}$, for each $n$.  It should be noted that we do not know the exact Puk\'anskzy invariant of these $\A_n$, only that $1$ and $2$ are members of $\puk{\A_n}$.
\end{eg}

We can also use Theorem \ref{Path.Cont.Result}, to show that the
Puk\'anszky invariant can not be used to give a continuous path of
non-conjugate singular masas even though the cardinality of the
set of non-conjugate singular masas is large enough. The proof is
omited, it can be found in \cite[Corollary 3.1.8]{Saw.Thesis}.
\begin{cor}
Let $\N$ be a separable \IIi factor. There is no continuous map
$t\mapsto\A(t)$ from $[0,1]$ into the set of all masas in $\N$
equiped with the $\dit$-metric such that $t\mapsto\puk{\A(t)}$ is
injective.
\end{cor}

\section{A $\Gamma$-invariant for masas}\label{Path.Gamma}
To show that all the uncountably many masas we shall produce are
pairwise non-conjugate via automorphisms of the underlying \IIi
factor, we introduce a conjugacy invariant.
\begin{dfn}\label{Path.Con.Invariant}
Let $\A$ be a masa in a \IIi factor $\N$.  Define $\Gamma(\A)$ to
be the supremum of $\tr(p)$ over all projections $p\in\A$ such
that $\A p$ contains non-trivial centralising sequences for $p\N
p$. If $\Gamma(\A)=0$, then we say that $\A$ is \emph{totally
non-$\Gamma$}.
\end{dfn}
Recall that a centralising sequence in a non-empty subset $B$ of a
\IIi factor $\N$ is a sequence $\{x_n\} \subset B$ with
\begin{align}
 \nm{x_ny -yx_n}_2 \rightarrow 0 \quad \text{for all} \quad y \in N . \notag
\end{align}
The centralising sequence $\{x_n\} \subset B$ is \emph{trivial} if
there is a sequence $\{\lambda_n \} \subset \mathbb{C}$ with
$\nm{x_n - \lambda_n }_2 \rightarrow 0$.

It is immediate that $\Gamma(\A)$ is a conjugacy invariant of
$\A$, in the sense that for an automorphism $\theta$ of $\N$, we
have $\Gamma(\theta(\A))=\Gamma(\A)$.

We shall produce masas in a similar fashion to Example
\ref{Path.Cont.Example}, taking a `direct sum' of a $\Gamma$-masa,
that is one containing non-trivial centralising sequences for its
underlying \IIi factor, and a totally non-$\Gamma$ masa.  The next
lemma is the tool that allows us to do this.
\begin{lem}\label{Path.Con.GammaSum}
Let $\A$ be a masa in a \IIi factor $\N$. Suppose that there is a
projection $p\in\A$ such that
\begin{itemize}
\item $\A p$ contains non-trivial centralising sequences for $p\N p$;
\item $\A(1-p)$ is totally non-$\Gamma$ in $(1-p)\N(1-p)$.
\end{itemize}
Then $\Gamma(\A)=\tr(p)$.
\begin{proof}
Take a projection $r\in\A$ such that $\A r$ contains non-trivial
centralising sequences for $r\N r$.  To obtain a contradiction,
suppose that $r\not\leq p$.  Let $(x_n)_{n=1}^\infty$ be a
non-trivial centralising sequence for $r\N r$ in $\A r$, write
$y_n=x_npr = x_nrp$ and $z_n=x_nr(1-p)$ so that $x_n=y_n+z_n$ for
all $n$. The sequence $(z_n)_{n=1}^\infty$ is a centralising
sequence of $r(1-p)\N r(1-p)$ and so is trivial by hypothesis.
Without losing generality, we may assume that $z_n=r(1-p)$ for all
$n$.

Take a partial isometry $v\in\N$ with $v^*v\leq r(1-p)$ and
$vv^*=p_0\leq pr$, so that $y_nv=x_nv$ and $v=vz_n=vx_n$.  Now
\begin{equation}\label{Puk.Uncount.GammaSum.1}
\nm{(y_n-1)p_0}_2=\nm{(y_n-1)v}_2=\nm{x_nv-vx_n}_2\rightarrow0,
\end{equation}
as $n\rightarrow\infty$.

By Kadison's Theorem on projections in a masa (\cite{Kadison.DiagMatrices}) choose orthogonal projections $(p_m)_{m=1}^{m_0}$ in $A$, with
$p_m\leq pr$ and $\tr(p_m)\leq\tr(r(1-p))$, for each $m$, so that
$\sum_{m=1}^{m_0}p_m=pr$.  Then, by
(\ref{Puk.Uncount.GammaSum.1}),
\[
 \nm{y_n-pr}_2 = \nm{(y_n - 1)pr}_2 \leq
 \sum_{m=1}^{m_0}\nm{(y_n-1)p_m}_2\rightarrow0, \notag
\]
so that $(x_n)_{n=1}^{\infty}$ is a trivial centralising sequence.
This contradiction ensures that $r\leq p$ and so
$\Gamma(\A)=\tr(p)$, as required.
\end{proof}
\end{lem}

The $\Gamma$-invariant is uniformally continuous with respect to
the $\dit$-metric on masas in separable \IIi factors.

\begin{lem}\label{Path.Con.GammaCts}
For masas $\A$ and $\B$ in a separable \IIi factor $\N$, we have
$$
\md{\Gamma(\A)-\Gamma(\B)}\leq 15\dit(\A,\B).
$$
\begin{proof}
Suppose that $\A$ and $\B$ are masas in $\N$ with
$\dit(\A,\B)\leq\epsilon$.  Let $u,p$ and $q$ be as in Theorem
\ref{Path.Cont.Pert}, so that
$$
\nm{1-p}_2=\nm{1-q}_2\leq15\epsilon.
$$
Given a projection $e\in\A$ such that $\A e$ has non-trivial
central sequences for $e\N e$, take $f=uepu^*$ --- a projection in
$\B q$ with $f \leq q$.  Since $u\A epu^*=\B f$, we can use
$u$ to conjugate the centralising sequences for $ep\R
ep$ lying in $\A e$ into centralising sequences for $f\R f$ lying in $\B f$.
Therefore,
$$
\Gamma(\B)\geq\tr(ep)=\tr(e)-\tr(e(1-p))\geq\tr(e)-\nm{e}_2
\nm{1-p}_2\geq \tr(e)-15\epsilon
$$
for every such projection $e \in \A$.  Hence,
$$
\Gamma(\B)\geq\Gamma(\A)-15\epsilon.
$$
By interchanging the roles of $\A$ and $\B$ we have
$$
\Gamma(\A)\geq\Gamma(\B)-15\epsilon,
$$
and these two inequalities combine to give the result.
\end{proof}
\end{lem}

One might attempt to produce uncountably many non-conjugate singular masas in the hyperfinite \IIi factor $\R$ with Puk\'anszky invariant $\{1\}$ by taking projections $e\in\R$ and singular masas $\B_1$ in $e\R e$ and $\B_2$ in $(1-e)\R(1-e)$ both with Puk\'anszky invariant $\{1\}$, such that $\B_1$ is $\Gamma$ in $e\R e$ and $\B_2$ is totally non-$\Gamma$ in $(1-e)\R(1-e)$.  The `direct-sum' $\A=\Set{b_1+b_2|b_1\in\B_1,b_2\in\B_2}$ will be a masa in $\R$ with $\Gamma(\A)=\tr(e)$ by Lemma \ref{Path.Con.GammaSum}.  Unfortunately, we do not have control over the exact Puk\'anszky invariant of such a masa $\A$, all we can say is that $1\in\puk{\A}$.  Indeed, there is a masa $\A$ in $\R$ with $\puk{\A}=\Set{1,2}$ for which there is a projection $e\in\A$ with $\tr(e)=1/2$ such that 
$$
\puk{\A e\subset e\R e}=\puk{\A(1-e)\subset (1-e)\R(1-e)}=\{1\}.
$$
Examples to this effect will be given in subsequent work by the second author.  In the next section, we get round this problem using Tauer masas to control the Puk\'anskzy invariant of these direct sums.
\section{A continuous path of singular masas}\label{Path.Con}
Here is the main result of this paper, from which Theorem
\ref{Path.Con.Cor} follows immediately.
\begin{thm}\label{Path.Con.Result}
There is a map $t\mapsto\A(t)$, taking each $t\in[0,1]$ to a masa
$\A(t)$ in $\R$ such that
\begin{enumerate}[(i)]
\item $\dit(\A(s),\A(t))\rightarrow 0$ as $\md{s-t}\rightarrow0$.\label{Path.Con.Result.1}
\item Every $\A(t)$ has Puk\'anszky invariant $\{1\}$.\label{Path.Con.Result.2}
\item Each $\A(t)$ is singular.\label{Path.Con.Result.3}
\item $\Gamma(\A(t))=t$, for each $t$.\label{Path.Con.Result.4}
\end{enumerate}
\end{thm}

We shall construct Tauer masas,
$\A(t)$, for a dense set of $t$ in $[0,1]$ with the required
properties, then use continuity to produce the required path. The
construction in the dense set of $t$ is based on a rapidly
increasing sequence of primes and adjusting the definition of the
approximately finite dimensional approximating algebras according
to $t$ being in suitable ranges of rationals.
\begin{notn}
Let $k_1=2$, and for each $r\geq 2$ take $k_r$ to be a prime
exceeding $k_1\dots k_{r-1}$.  Let $\M_r$ to be the algebra  of
$k_r\times k_r$ matrices.  By \cite[Theorem 3.2]{Popa.Orth}, there
is a family $(^r\D^{(m)})_{m=0}^{k_1\dots k_{r-1}}$ of pairwise
orthogonal masas in $\M_r$.  Write $^re^{(m)}_l$ for the minimal
projections of $^r\D^{(m)}$ indexed by $l=0,1\dots,k_r-1$. Let
$\N_n$ be the tensor product $\bigotimes_{r=1}^n\M_r$.  We have
the natural unital inclusion $x\mapsto x\otimes 1$ of $\N_n$
inside $\N_{n+1}$ and we work in the hyperfinite \IIi factor $\R$,
obtained as the direct limit of these $\N_n$ with respect to
normalised trace.

For each $n\in\mathbb N$ write
$$
I_n=\Set{\frac{m}{k_1\dots k_n}|m=0,1,2,\dots,k_1\dots k_n},
$$
so that $I_n\subset I_{n+1}$, for each $n$.  Let
$I=\bigcup_{n=1}^\infty I_n$ --- a dense set of rationals in $[0,1]$.
For each $t\in I$, we will define a Tauer masa $\A(t)$ in $\R$ with respect to the chain $(\N_n)_{n=n_0(t)}^\infty$, where $n_0(t)$ is the minimal $n$ for which $t\in I_n$.  For each $n\geq n_0(t)$, we denote the $n$-th approximate of $\A(t)$ by $\A_n(t)$, and enumerate the minimal projections of $\A_n(t)$ as $^nf_m(t)$ for $0\leq m<k_1\dots k_n$.
\end{notn}

\begin{con}\label{Path.Con.Dense}
The process begins by defining
$\A_0(0)=\A_0(1/2)=\A_0(1)={}^1\D^{(0)}$ with the minimal
projections $^1f_m(0)={}^1f_m(1/2)={}^1f_m(1)={}^1e^{(0)}_m$
coinciding for $m=0,1$.  For some $n_1$, suppose that we have
defined $\A_n(t)$ and enumerated the minimal projections $^nf_m(t)$, for all
$t\in I_{n_1}$ and $n_0(t)\leq n\leq n_1$.  For $t\in I_{n_1}$,
the definition of $\A_{n_1+1}(t)$ is split into two cases,
depending on whether $n_1$ is even or odd.

\begin{enumerate}
\item $n_1$ is \emph{even}: Set
\begin{equation}\label{Puk.Uncount.Dense.9}
\A_{n_1+1}(t)=\bigoplus_{m=0}^{k_1\dots k_{n_1}-1}{}^{n_1}f_m(t)\otimes{}^{n_1+1}\D^{(m)}.
\end{equation}
Enumerate the minimal projections $^{n_1+1}f_{m'}(t)$ by dividing $m'$ by $k_{n_1+1}$ to obtain $m'=k_{n_1+1}m+l$ for some $0\leq
l<k_{n_1+1}$.  Now take
\begin{equation}\label{Puk.Uncount.Dense.3}
^{n_1+1}f_{m'}(t)={}^{n_1}f_m(t)\otimes{}^{n_1+1}e^{(m)}_l.
\end{equation}
\item $n_1$ is \emph{odd}: Here we take
\begin{align}
\A_{n_1+1}(t)=&\bigoplus_{m=0}^{tk_1\dots k_{n_1}-1}{}^{n_1}f_m(t)
\otimes{}^{n_1+1}\D^{(k_1\dots k_{n_1})}\nonumber\\
&\quad\oplus\bigoplus_{m=tk_1\dots k_{n_1}}^{k_1\dots k_{n_1}-1}{}^{n_1}f_m(t)\otimes{}^{n_1+1}\D^{(m)}.\label{Puk.Uncount.Dense.8}
\end{align}
The enumeration of the minimal projections happens in the same way
as the even $n_1$ case.  Namely, given $0\leq m'<k_1\dots k_{n_1+1}$ write
$m'=mk_{n_1+1}+l$ for some $0\leq l<k_{n_1+1}$ and set
\begin{equation}\label{Puk.Uncount.Dense.4}
^{n_1+1}f_{m'}(t)=\left\{\begin{array}{ll}^{n_1}f_m(t)\otimes{}^{n_1+1}e^{(k_1\dots k_{n_1})}_l&0\leq m<tk_1\dots k_{n_1}\\^{n_1}f_m(t)\otimes{}^{n_1+1}e^{(m)}_l&tk_1\dots k_{n_1}\leq m<k_1\dots k_{n_1}\end{array}\right..
\end{equation}
\end{enumerate}

It remains to define $\A_{n_1+1}(t)$ when $t\in I_{n_1+1}\setminus
I_{n_1}$.  In this case this is the first approximate of the Tauer masa $\A(t)$.  Write
$m_0=\floor{tk_1\dots k_{n_1}}$ and define the minimal projections of $\A_{n_1+1}(t)$ by
\begin{equation}\label{Path.Con.StartStage}
^{n_1+1}f_m(t)=\left\{\begin{array}{ll}^{n_1+1}f_m((m_0+1)/k_1\dots k_{n_1})&0\leq m<tk_1\dots k_{n_1+1}\\^{n_1+1}f_m(m_0/k_1\dots k_{n_1})&tk_1\dots k_{n_1+1}\leq m<k_1\dots k_{n_1+1}\end{array}\right..
\end{equation}
\end{con}

Theorem 4.1 of \cite{Saw.Tauer} shows that the Tauer masas
constructed above have $\puk{\A(t)}=\{1\}$, which is
condition (\ref{Path.Con.Result.2}) of Theorem \ref{Path.Con.Result}.  We now check that these masas satisfy conditions (\ref{Path.Con.Result.3}) and (\ref{Path.Con.Result.4}) of Theorem \ref{Path.Con.Result}.
\begin{lem}\label{Path.Con.Sing}
The Tauer masas $\A(t)$ of Construction \ref{Path.Con.Dense} are singular.
\begin{proof}
Fix $t\in I$ and let $n\geq n_0(t)$ be even.  In the notation of
(\ref{Path.Intro.Inclusion}), the even stage
of Construction \ref{Path.Con.Dense} gives
$$
\A_{n+1,n}^{\left({}^nf_m(t)\right)}(t)={}^{n+1}\D^{(m)}.
$$
Take a unitary $w\in{}^{n+1}\D^{(m)}$ with $\tr(w)=0$.  When $m'\neq
m$, the orthogonality of $^{n+1}\D^{(m)}$ and $^{n+1}\D^{(m')}$ gives $\ce{^{n+1}\D^{(m')}}{w}=0$.  The singularity of $\A(t)$ then follows from Proposition \ref{Path.Intro.Sing}.
\end{proof}
\end{lem}

The next Lemma verifies the hypothesis of Lemma \ref{Path.Con.GammaSum}, so the masas of Construction \ref{Path.Con.Dense} have $\Gamma(\A(t))=t$.
\begin{lem}\label{Path.Con.Gamma}
Fix $t\in I$ and write $n_0$ for $n_0(t)$.  Let
\begin{equation}\label{Puk.Uncount.Defp}
p=\sum_{m=0}^{tk_1\dots k_{n_0}-1}{}^{n_0}f_m(t),
\end{equation}
a projection in $\A(t)$.  Then
\begin{enumerate}
\item $\A(t)p$ contains non-trivial centralising sequences for $p\R p$;
\item $\A(t)(1-p)$ is totally non $\Gamma$ in $(1-p)\R(1-p)$.
\end{enumerate}
\end{lem}

\begin{proof}[Proof of 1:]
Note that
$$
p=\sum_{m=0}^{tk_1\dots k_n-1}{}^{n}f_m(t),
$$
for all $n\geq n_0$.  Fix $n\geq n_0$ odd and consider $x_1,\dots
,x_r\in\N_n$.  Let $v\in{}^{n+1}\D^{(k_1\dots k_n)}$ be a unitary
with $\tr(v)=0$.  Examining the odd $n$ form of Construction
\ref{Path.Con.Dense}, we see that
$$
u=\sum_{m=0}^{tk_1\dots k_n-1}{}^{n}f_m(t)\otimes v=p\otimes
v\in\N_n\otimes\M_{n+1}=\N_{n+1}
$$
is a trace free unitary in $\A_{n+1}(t)p$.  It is then immediate
that $u$ commutes with each $px_ip$, and so $\A(t)p$ contains
non-trivial centralising sequences for $p\R p$ by the $\nm{.}_2$-density of
$\cup_{n=1}^\infty\N_n$ in $\R$.
\end{proof}

We prove part 2 of Lemma \ref{Path.Con.Gamma} in two stages.
We first establish an orthogonality condition which suffices to
establish that no $\A e$ can contain centralising sequences for
$e\R e$, when $e\leq 1-p$ is a minimal projection of some
$\A_n(t)$. A density argument, which contains the proof of an
observation of Popa (\cite[Remark 5.4.2]{Popa.Orth}, also found in \cite[Lemma
2.1]{Bisch.CentralSeq2}), then completes the proof of Lemma
\ref{Path.Con.Gamma}.
\begin{lem}\label{Puk.Uncount.Orth}
Fix $t\in I$, $n\geq n_0(t)$ and $m,m'$ with $tk_1\dots k_n\leq
m<m'<k_1\dots k_n$.  Let $v$ be a partial isometry in $\N_n$ with
$vv^*=\sp^nf_m(t)$ and $v^*v=\sp^nf_{m'}(t)$. Then
$v(\A(t)\sp^nf_{m'}(t))v^*$ is orthogonal to $\A(t)\sp^nf_m(t)$ in
$^nf_m(t)\R\sp^nf_m(t)$.
\begin{proof}
Fix $n\geq n_0$ and regard $\R$ as $\N_n\otimes\R_1$, where $\R_1$
is generated as the infinite von Neumann tensor product
$(\bigotimes_{r=n+1}^\infty\M_r)''$ with respect to the unique
normalised trace.  Using the notation of (\ref{Path.Intro.Inclusion}), for
$n_1>n$ we have
$$
\A_{n_1}(t)=\bigoplus_{m=0}^{k_1\dots k_n-1}{}^nf_m(t)
\otimes\A_{n_1,n}^{\left({}^nf_m(t)\right)}(t),
$$
for masas $\A_{n_1,n}^{\left({}^nf_m(t)\right)}(t)$ in
$\bigotimes_{r=n+1}^{n_1}\M_r$.  In this way we obtain Tauer masas
$$
\A_{\infty,n}^{\left({}^nf_m(t)\right)}(t)=
\left(\bigcup_{n_1=n+1}^\infty\A_{n_1,n}^{\left({}^nf_m(t)\right)}(t)\right)''
$$
in $\R_1$, so that
$$
\A(t)=\bigoplus_{m=0}^{k_1\dots k_n-1}{}^nf_m(t)
\otimes\A_{\infty,n}^{\left({}^nf_m(t)\right)}(t).
$$
Now take $m,m'$ and $v$ as in the statement, and note that
$$
v\A(t){}^nf_{m'}(t)v^*={}^nf_m(t)\otimes\A_{\infty,n}^{\left({}^nf_{m'}(t)\right)}(t),
$$
so that it suffices to show that
$\A_{\infty,n}^{\left({}^nf_{m}(t)\right)}(t)$ and
$\A_{\infty,n}^{\left({}^nf_{m'}(t)\right)}(t)$ are orthogonal masas
in $\R_1$.  We shall show that
$\A_{n_1,n}^{\left({}^nf_{m}(t)\right)}(t)$ and
$\A_{n_1,n}^{\left({}^nf_{m'}(t)\right)}(t)$ are orthogonal in
$\bigotimes_{r=n+1}^{n_1}\M_r$, for all $n_1>n$, from which the result immediately
follows by density.

To this end, note that Construction \ref{Path.Con.Dense} gives
$\A_{n+1,n}^{\left({}^nf_{m}(t)\right)}(t)={}^{n+1}\D^{(m)}$ and
$\A_{n+1,n}^{\left({}^nf_{m'}(t)\right)}(t)={}^{n+1}\D^{(m')}$, from
(\ref{Puk.Uncount.Dense.9}) when $n$ is even and from
(\ref{Puk.Uncount.Dense.8}) when $n$ is odd.  In the latter case, we use the hypothesis that $tk_1\dots k_n\leq m<m'$.  As
$\D^{(m)}$ and $\D^{(m')}$ are orthogonal masas in $\M_{n+1}$, the claim holds when $n_1=n+1$.

Suppose inductively that the claim holds for some $n_1>n$.  Write
$$
\A_{n_1+1,n}^{\left({}^nf_m(t)\right)}(t)
=\bigoplus_{g\in\Proj\left(\A_{n_1,n}^{\left({}^nf_m(t)\right)}(t)\right)}g\otimes\B^{(g,m)},
$$
and
$$
\A_{n_1+1,n}^{\left({}^nf_{m'}(t)\right)}(t)
=\bigoplus_{h\in\Proj\left(\A_{n_1,n}^{\left({}^nf_{m'}(t)\right)}(t)
\right)}h\otimes\B^{(h,m')},
$$
for masas $B^{(g,m)}$ and $B^{(h,m')}$ in $\M_{n_1+1}$.  Again,
Construction \ref{Path.Con.Dense} ensures that all these masas
are pairwise orthogonal.  This is immediate from
(\ref{Puk.Uncount.Dense.9}) for even $n_1$; when $n_1$ is odd we
again use the hypothesis $tk_1\dots k_n\leq m<m'$ in our examination
of (\ref{Puk.Uncount.Dense.8}).  The orthogonality of
$\A_{n_1+1,n}^{\left({}^nf_m(t)\right)}(t)$ and
$\A_{n_1+1,n}^{\left({}^nf_{m'}(t)\right)}(t)$ follows immediately,
yielding the result.
\end{proof}
\end{lem}

\begin{proof}[Proof of part 2 of Lemma \ref{Path.Con.Gamma}:]
Take $t\in I$ and fix some projection $0\neq e\leq 1-p$ in $\A(t)$.
For each $n\in\mathbb N$, find $l_n\geq n_0(t)$ and a family
$P_n\subset\Proj(\A_{l_n}(t))$ of minimal projections in
$\A_{l_n}(t)$ lying under $1-p$, such that upon writing $q_n=\sum_{q\in
P_n}q$, we have
$$
\nm{q_n-e}_2^2<1/n.
$$

For each $n$, take a permutation $\sigma_n$ of $P_n$ with no fixed points.  Take partial isometries $v_{\sigma_n(q),q}$ in $\N_{l_n}$ with $v_{\sigma_n(q),q}{v_{\sigma_n(q),q}}^*=\sigma_n(q)$ and ${v_{\sigma_n(q),q}}^*v_{\sigma_n(q),q}=q$.  Define 
$$
x_n=\sum_{q\in P_n}v_{\sigma_n(q),q}+(1-q_n),
$$
a unitary in $\N_{l_n}$ which has $x_nqx_n^*=\sigma_n(q)$, for every $q\in P_n$.  Observe that
$$
x_n(\A q_n)x_n^*=\bigoplus_{q\in P_n}x_n(\A q)x_n^*=\bigoplus_{q\in P_n}v_{\sigma_n(q),q}(\A q){v_{\sigma_n(q),q}}^*=\bigoplus_{q\in P_n}v_{\sigma_n(q),q}\A {v_{\sigma_n(q),q}}^*
$$
is orthogonal to $\bigoplus_{q\in P_n}\A\sigma_n(q)=\A q_n$ in $q_n\R q_n$ by Lemma \ref{Puk.Uncount.Orth}.

Suppose that $\A e$ contains non-trivial centralising sequences for
$e\R e$.  Find a sequence of unitaries $u_n\in\A$, with
$\tr(u_ne)=0$ for each $n$, and such that
\begin{equation}\label{Puk.Uncount.Proof2.1}
\nm{eu_nex_ne-ex_neu_ne}_2<\nm{e-q_n}_2.
\end{equation}
We have the following simple estimate, showing that $u_nq_n$
asymptotically commutes with the $q_nx_nq_n$:
\begin{align*}
&\nm{q_nu_nq_nx_nq_n-q_nx_nq_nu_nq_n}_2\\
\leq&\nm{(q_n-e)u_nq_nx_nq_n}_2+\nm{eu_n(q_n-e)x_nq_n}_2+\nm{eu_nex_n(q_n-e)}_2\\
&+\nm{eu_nex_ne-ex_neu_ne}_2+\nm{ex_neu_n(e-q_n)}_2+\nm{ex_n(e-q_n)u_nq_nm}_2\\
&+\nm{(e-q_n)x_nq_nu_nq_n}_2\\
\leq&\sp7\nm{e-q_n}_2\rightarrow0.
\end{align*}
On the other hand, using $x_nq_n=q_nx_n$ we have
\begin{align*}
&\sp\nm{q_nx_nq_nu_nq_n-q_nu_nq_nx_nq_n}^2_2\\
=&\sp\nm{q_nx_nu_nq_nx_n^*q_n-u_nq_n}_2^2\\
=&\sp\nm{q_nx_nu_nq_nx_n^*q_n}_2^2+\nm{u_nq_n}_2^2-2\Re\tr(x_nu_nq_nx_n^*u_nq_n)\\
=&\sp2\nm{q_n}_2^2-2\Re\tr(x_nu_nq_nx_n^*)\tr(u_nq_n)/\tr(q_n)\rightarrow2\nm{e}_2^2\neq 0,
\end{align*}
where the last line comes from the orthogonality of $x_n(\A q_n)x_n^*$ and $\A q_n$ in $q_n\R q_n$ --- the quotient of
$\tr(q_n)$ appearing as a normalisation constant.  The convergence is a
simple calculation, as
$$
\md{\tr(u_nq_n)}\leq\md{\tr(u_ne)}+\md{\tr(u_n(q_n-e))}\leq 0+\nm{u_n}_2\nm{q_n-e}_2\rightarrow 0.
$$
This contradiction completes the proof.
\end{proof}

For $t$ in the dense subset $I$ of $[0,1]$, we have singular Tauer masas $\A(t)$ with $\Gamma(\A(t))=t$.  We wish to use completeness
to define $\A(t)$ for $t\in[0,1]\setminus I$ and so we need to
control the distance between the $\A(t)$'s we have already defined.
It is here that the form of $\A_{n_0(t)}(t)$ specified in
Construction \ref{Path.Con.Dense} becomes relevant.
\begin{lem}\label{Path.Con.Same}
Fix $s,t\in I$ with $s<t$.  Let $n_0$ be the maximum of $n_0(s)$ and $n_0(t)$ and take
$$
q=\sum_{m=0}^{sk_1\dots k_{n_0}-1}{}^{n_0}f_m(s)+\sum_{m=tk_1\dots k_{n_0}}^{k_1\dots k_{n_0}-1}{}^{n_0}f_m(s)
$$
a projection of trace $1-(t-s)$.  This $q$ lies in $\A(s)\cap\A(t)$ and $\A(s)q=\A(t)q$.
\begin{proof}
We shall demonstrate that Construction \ref{Path.Con.Dense}
ensures that whenever we have $s,t\in I_n$, then
\begin{equation}\label{Puk.Uncount.Dense.1}
^nf_m(s)={}^nf_m(t),
\end{equation}
for all $m$ with
\begin{equation}\label{Puk.Uncount.Dense.2}
0\leq m<sk_1\dots k_n\textrm{ or }tk_1\dots k_n\leq m<k_1\dots k_n.
\end{equation}
This will immediately show that $q$ lies in $\A(t)$, as well as
$\A(s)$.  Furthermore, as $\A(s)q$ and $\A(t)q$ are generated by all $^nf_m(s)$ and $^nf_m(t)$ respectively, with
$n\geq\max\{n_0(s),n_0(t)\}$ and $m$ satisfying
(\ref{Puk.Uncount.Dense.2}), the claim also implies that
$\A(s)q=\A(t)q$, as required.

We proceed by induction on $n$.  When $n=1$, the
result is certainly true, as Construction \ref{Path.Con.Dense}
began by defining $^1f_m(0)=\sp^1f_m(1/2)=\sp^1f_m(1)$ for $m=0,1$.  Suppose that we have established the
claim for all $n\leq n_1$.  We investigate the $n_1+1$ situation,
starting with the case when $s$ and $t$ both lie in $I_{n_1}$.

Take $s,t\in I_{n_1}$ with $s<t$.  Take $m'$ with either $0\leq
m'<sk_1\dots k_{n_1+1}$ or $tk_1\dots k_{n_1+1}\leq m'<k_1\dots
k_{n_1+1}$, and divide by $k_{n_1+1}$ to obtain $m'=mk_{n_1+1}+l$
with $0\leq l<k_{n_1+1}$.  This $m$ must have $0\leq m<sk_1\dots
k_{n_1}$ in the first case or $tk_1\dots k_{n_1}\leq m<k_1\dots
k_{n_1}$ in the second.  In any event, the inductive hypothesis
ensures that $^{n_1}f_m(s)={}^{n_1}f_m(t)$.  When $n_1$ is even, the
definition (\ref{Puk.Uncount.Dense.3}) of $^{n_1+1}f_{m'}(s)$ and
$^{n_1+1}f_{m'}(t)$ immediately gives
$^{(n_1+1)}f_{m'}(s)={}^{n_1+1}f_{m'}(t)$.  When $n_1$ is odd, this
is also true, as we have excluded the possibility that $sk_1\dots
k_{n_1}\leq m<tk_1\dots k_{n_1}$, so both these minimal projections
must come from the same case of equation
(\ref{Puk.Uncount.Dense.4}).  Therefore, the minimal
projections $^{n_1+1}f_{m'}(s)$ and $^{n_1+1}f_{m'}(t)$ coincide
whenever they are required to do so.

We now examine what happens when precisely one of $s$ and $t$ lies in
$I_{n_1+1}\setminus I_{n_1}$.  Take $s$ in $I_{n_1}$ and $t\in I_{n_1+1}\setminus
I_{n_1}$ with $s<t$.  As in the definition of $\A_{n_1+1}(t)$, we write
$m_0=\floor{tk_1\dots k_n}$ so that $s\leq m_0/k_1\dots k_{n_1}$.
For $0\leq m< sk_1\dots k_{n_1+1}$, we have
$$
^{n_1+1}f_m(s)={}^{n_1+1}f_m((m_0+1)/k_1\dots k_{n_1})={}^{n_1+1}f_m(t),
$$
where the second equality is the definition, (\ref{Path.Con.StartStage}), of $^{n_1+1}f_m(t)$, and the first follows as the $m$-th minimal
projections for $\A_{n_1+1}(s)$ and $\A_{n_1+1}((m_0+1)k_1\dots
k_{n_1})$ coincide by the case we analysed in the previous paragraph.  When
$tk_1\dots k_{n_1+1}\leq m<k_1\dots k_{n_1+1}$, we have
$$
^{n_1+1}f_m(t)={}^{n_1+1}f_m(m_0/k_1\dots k_{n_1})={}^{n_1+1}f_m(s),
$$
the first equality being (\ref{Path.Con.StartStage}) -- the definition of $^{n_1+1}f_m(t)$, and the
second equality is (\ref{Puk.Uncount.Dense.1}) for appropriate
minimal projections of $\A_{n_1+1}(s)$ and $\A_{n_1+1}(m_0/k_1\dots
k_{n_1})$, as $m\geq m_0k_{n_1+1}$. These last two algebras may
turn out to be the same, but then the minimal projections will
certainly coincide.  Interchanging the roles of $s$ and $t$ above
ensures that the claim holds for $n_1+1$ whenever either $s$ or $t$ lies in $I_{n_1}$.

We complete the proof by examining the situation when $s,t\in
I_{n_1+1}\setminus I_{n_1}$.  Take $s<t$
with $s,t\in I_{n_1+1}\setminus I_{n_1}$.  Suppose first that
$\floor{sk_1\dots k_{n_1}}=\floor{tk_1\dots k_{n_1}}=m_0$.  In this instance the definition, (\ref{Path.Con.StartStage}), of the minimal projections $^{n_1+1}f_m(s)$ and $^{n_1+1}f_m(t)$ ensures that these projections conincide for all $m$ with $0\leq m<sk_1\dots k_{n_1+1}$ or $tk_1\dots k_{n_1+1}\leq m<k_1\dots k_{n_1+1}$.

Finally, suppose that $\floor{sk_1\dots k_{n_1}}=m_0<
m_1=\floor{tk_1\dots k_{n_1}}$.  Given $m$ with $0\leq m<sk_1\dots k_{n_1+1}$, (\ref{Path.Con.StartStage}) ensures that $^{n_1+1}f_m(s)=\sp^{n_1+1}f_m(m_0+1)$ and $^{n_1+1}f_m(t)=\sp^{n_1+1}f_m(m_1)$.  Since $m<sk_1\dots k_{n_1+1}<(m_0+1)k_1\dots k_{n_1+1}$, the case when $s,t\in I_n$ (with $s=m_0+1$ and $t=m_1$) ensures that $^{n_1+1}f_m(m_0+1)=\sp^{n_1+1}f_m(m_1)$.  In conclusion, we have $^{n_1+1}f_m(s)=\sp^{n_1+1}f_m(t)$ as required.  The case when $tk_1\dots k_{n_1+1}\leq m<k_1\dots k_{n_1+1}$ is similar, and this completes the proof.
\end{proof}
\end{lem}

\begin{cor}\label{Path.Con.Dist}
For $s,t\in I$ we have
$$
\nmit{\Ce{\A(s)}-\Ce{\A(t)}}\leq 2\sqrt{\md{s-t}}
$$
\begin{proof}
We may assume that $s<t$.  Let $n_0$ be the maximum of $n_0(s)$ and
$n_0(t)$.  Let $q$ be the projection of Proposition
\ref{Path.Con.Same}, so that $\A(s)q=\A(t)q$.   The simple estimate
$$
\nmit{\Ce{\A(s)}-\Ce{\A(t)}}\leq 2\nm{1-q}_2,
$$
can be found in part (i) of Theorem 6.5 in \cite{Sinclair.PertSubalg}.
As $\tr(1-q)=t-s$, this is exactly what was claimed.
\end{proof}
\end{cor}

We can now combine the results of this section to esablish Theorem
\ref{Path.Con.Result}.
\begin{proof}[Proof of Theorem \ref{Path.Con.Result}.]
For $t\in I$ we take $\A(t)$ to be the Tauer masa produced in Construction \ref{Path.Con.Dense}.  When $t\in[0,1]\setminus I$, we
define $\A(t)$ by taking a sequence $t_n\rightarrow t$ with each
$t_n$ in the dense set of rationals $I$.  The resulting sequence of
masas $(\A(t_n))_{n=1}^\infty$ is then $\dit$-Cauchy by Corollary
\ref{Path.Con.Dist}, and so converges to a masa $\A(t)$ in $\R$.  Recall that the set of all von Neumann subalgebras of $\R$ is a $\dit$-complete metric space, and the subset of masas is closed, \cite{Christensen.Subalgebras}.  This masa is well
defined, in that $\A(t)$ is independent of the choice of sequence
$t_n$ in $I$ converging to $t$.  An approximation argument extends Corollary \ref{Path.Con.Dist} to show that $\dit(\A_s,\A_t)\rightarrow 0$ whenever $\md{s-t}\rightarrow 0$.

Furthermore each $A(t)$ is singular, as this holds for $t\in I$ (Lemma \ref{Path.Con.Sing}) and
the set of all singular masas is closed; again this can be found in \cite{Christensen.Subalgebras}. All the $\A(t)$ have Puk\'anszky invariant $\{1\}$; for $t\in I$ this is Theorem 4.1 of \cite{Saw.Tauer} and Corollary \ref{Path.Cont.Closed} then gives the result for general $t$. That $\Gamma(\A(t))=t$ for
every $t\in[0,1]$ follows first by observing that Lemma
\ref{Path.Con.Gamma} combines with Lemma \ref{Path.Con.GammaSum} to
give the result for $t\in I$. Continuity gives the result for all
$t$, this time in the form of Lemma
\ref{Path.Con.GammaCts}.
\end{proof}

\end{document}